\documentclass{fic-l}

\newtheorem{theorem}{Theorem}

\newtheorem{corollary}[theorem]{Corollary}

\newtheorem{fact}{Fact}

\theoremstyle{definition}

\theoremstyle{remark}
\newtheorem*{remark*}{Remark}

\numberwithin{equation}{section}

\begin{document}

\title{On sequentially $h$-complete groups}

\author{G\'{a}bor Luk\'{a}cs}
\address{Department of Mathematics \& Statistics, York University\\ 
4700 Keele Street, Toronto, Ontario, M3J 1P3, Canada}
\curraddr{FB 3 - Mathematik und Informatik, Universit\"{a}t Bremen\\ 
Bibliothekstrasse 1, 28359 Bremen, Germany}
\email{lukacs@mathstat.yorku.ca}
\thanks{I gratefully acknowledge the financial support received from York 
University that enabled me to do this research.}


\begin{abstract}
A topological group $G$ is {\em sequentially $h$-complete} if all the
continuous homomorphic images of $G$ are sequentially complete. In this
paper we give necessary and sufficient conditions on a complete group for
being compact, using the language of sequential $h$-completeness. In the
process of obtaining such conditions, we establish a structure theorem for
$\omega$-precompact sequentially $h$-complete groups. As a consequence we
obtain a reduction theorem for the problem of $c$-compactness.
\end{abstract}

\maketitle

\nocite{Pontr}

All topological groups in this paper are assumed to be Hausdorff.

A topological group $G$ is {\em sequentially $h$-complete} if all the
continuous homomorphic images of $G$ are sequentially complete (i.e.,
every Cauchy-sequence converges). $G$ is called {\em precompact} if for
any neighborhood $U$ of the identity element there exists a finite subset
$F$ of $G$ such that $G=UF$.

In \cite[Theorem 3.6]{DikTka} Dikranjan and Tka{\v{c}}enko proved that
nilpotent sequentially $h$-complete groups are precompact (also see
\cite{Dik2}).  Thus, if a group is nilpotent, sequentially $h$-complete
and complete, then it is compact.

Inspired by this result, the aim of this paper is to give necessary and 
sufficient conditions on a complete group for being compact, using the 
language of sequential $h$-completeness. This aim is carried out in 
Theorem~\ref{thm:main}. 

For an infinite cardinal $\tau$, a topological group $G$ is 
{\em $\tau$-precompact} if for any neighborhood $U$ of the 
identity element there exists $F \subset G$ such that $G = UF$ and 
$|F| \leq \tau$. In order to prove Theorem~\ref{thm:main}, we will first 
establish a strengthened version of the Guran's Embedding Theorem 
for $\omega$-precompact sequentially $h$-completely groups 
(Theorem~\ref{thm:struct}).

A topological group $G$ is {\em $c$-compact} if for any topological group
$H$ the projection $\pi_H:G \times H \rightarrow H$ maps closed subgroups
of $G \times H$ onto closed subgroups of $H$ (see \cite{Manes}, \cite{DikGiu} 
and \cite{CGT}, as well as \cite{CleTho}). The problem of whether every
$c$-compact topological group is compact has been an open question for
more than ten years. As a consequence of Theorem~\ref{thm:main}, we obtain
that the problem of $c$-compactness can be reduced to the second-countable
case (Theorem~\ref{thm:reduction}).

\bigskip

The following Theorem is a slight generalization of Theorem~3.2 from
\cite{DikUsp}:

\begin{theorem} \label{thm:open}
Let $G$ be an $\omega$-precompact sequentially $h$-complete topological 
group. Then every continuous homomorphism $f: G \rightarrow H$ onto a 
group $H$ of countable pseudocharacter is open. 
\end{theorem}

In order to prove Theorem~\ref{thm:open}, we need the following three
facts, two of which are due to Guran:

\begin{fact}[Guran's Embedding Theorem] \label{fact:Guran} A topological 
group is $\tau$-pre\-com\-pact if and only if it is topologically 
isomorphic to a subgroup of a direct product of topological groups of weight 
$\le\tau$. {\rm (Theorem~4.1.3 in \cite{TkaTop1}.)}
\end{fact}

\begin{fact}[Banach's Open Map Theorem] \label{fact:Banach}
Any continuous homomorphism from a separable complete metrizable group 
onto a Baire group is open. {\rm (Corollary V.4 in \cite{Hus}.)} 
\end{fact}

\begin{fact} \label{fact:Guran2}
Let $G$ be an $\omega$-precompact topological group of countable 
pseudocharacter. Then $G$ admits a coarser second countable group 
topology. {\rm (Corollary 4 in \cite{Guran}.)}
\end{fact}

The proof below is just a slight modification of the proof of 
Theorem~3.2 from \cite{DikUsp} mentioned above:

\begin{proof}
First, suppose that $H$ is metrizable.
Let $U$ be a neighborhood of $e$ in $G$. Since, by 
Fact~\ref{fact:Guran}, $G$ embeds into a product of separable metrizable 
group, we may assume that $U  = g^{-1}(V)$ for some continuous homomorphism 
$g: G \rightarrow M$ onto a separable metrizable group $M$ and some 
neighborhood $V$ of $e$ in $M$. 
Let $h=(f,g): G \rightarrow H \times M$, and put $L=h(G)$. Let 
$p: L \rightarrow H$ and $q: L \rightarrow  M$ be the restrictions of the 
canonical projections $H \times M \rightarrow H$ and 
$H \times M \rightarrow M$. Clearly, one has $f=ph$ and $g=qh$. Since $q$ 
is continuous, $W = q^{-1}(V)$ is open. We have 
$h(U) = h (g^{-1}(V)) =h (h^{-1} q^{-1} (V)) \nolinebreak=\nolinebreak W$ 
and $f(U) = ph(U) = p(W)$.
Since $G$ is sequentially $h$-complete, the groups $L$ and $H$, being 
homomorphic images of $G$, are sequentially complete. Since $H$ and 
$L$ are metrizable, this means that they are simply complete. They are 
also separable, because they are metrizable and $\omega$-precompact. Thus, 
by Fact~\ref{fact:Banach}, $p: L \rightarrow H$ is open, and hence
$f(U) = p(W)$ is open in $H$.

To show the general case, suppose that the topology $\mathcal{T}$ of the 
group $H$ is of countable pseudocharacter;  then $H$ admits a coarser  
second countable topology  $\mathcal{T}^\prime$ (by Fact~\ref{fact:Guran2}); 
put $\iota: (H,\mathcal{T}) \rightarrow (H,\mathcal{T}^\prime)$ to be the 
identity map. Since $\mathcal{T}^\prime$ is metrizable, the continuous 
homomorphism $\iota \circ \varphi: G \rightarrow (H,\mathcal{T}^\prime)$ 
is open by what we have already proved, and thus $\varphi$ is also open.
\end{proof}

\begin{corollary} \label{cor:open:countpseudo}
Let $G$ be an $\omega$-precompact  sequentially $h$-complete group. The 
following statements are equivalent:

\begin{list}{{\rm (\roman{enumi})}}
{\usecounter{enumi}\setlength{\labelwidth}{25pt}\setlength{\topsep}{2pt}
\setlength{\itemsep}{0pt} \setlength{\leftmargin}{20pt}}

\item
$G$ is second countable;

\item
$G$ contains a countable network;

\item 
$G$ has a countable pseudocharacter;

\item $G$ is metrizable.

\end{list}
\end{corollary}

\begin{proof}
The implications (i) $\Rightarrow$ (ii) $\Rightarrow$ (iii) are 
obvious, and so is the equivalence of (i) and (iv).

(iii) $\Rightarrow$ (iv):  If $(G,\mathcal{T})$ is of countable 
pseudocharacter, it admits a coarser  second countable topology 
$\mathcal{T}^\prime $(by Fact~\ref{fact:Guran2}); put 
$\iota: G \rightarrow G_1$ to be the identity map. By 
Theorem~\ref{thm:open}, $\iota$ is open, and thus 
$\mathcal{T}$ is second countable, as desired.
\end{proof}

\bigskip

A topological group $G$ is {\em minimal} if every continuous isomorphism
$\varphi : G \rightarrow H$ is a homeomorphism, or equivalently, if the
topology of $G$ is a coarsest (Hausdorff) group topology on $G$. 
The group $G$ is {\em totally minimal} if every continuous surjective
homomorphism $f: G \rightarrow H$ is open; in other words, $G$ is totally
minimal if every quotient of $G$ is minimal. 

\begin{corollary} \label{cor:totmin-gen}
Every $\omega$-precompact sequentially $h$-complete topological group of 
countable pseudocharacter is totally minimal and metrizable.
\end{corollary}

\nocite{Arh}

\begin{proof}
Let $G$ be an $\omega$-precompact sequentially $h$-complete group of 
countable pseudocharacter. By Corollary~\ref{cor:open:countpseudo}, $G$ is 
metrizable and contains a countable network. Sequential $h$-completeness 
and the property of having a countable network are preserved under 
continuous homomorphic images, so for every  continuous 
homomorphism $\varphi:\nolinebreak G \rightarrow H$ onto a topological 
group $H$, the group $H$ is sequentially $h$-complete and contains a 
countable network; in particular, $H$ is $\omega$-precompact. Thus, by 
Corollary~\ref{cor:open:countpseudo}, $H$ is of countable 
pseudocharacter. Therefore, by Theorem~\ref{thm:open}, $\varphi$ is open.
\end{proof}

Corollary~\ref{cor:totmin} generalizes \cite[3.3]{DikUsp} to the
sequentially $h$-complete topological groups.

\begin{corollary} \label{cor:totmin}
Every sequentially $h$-complete topological group with a countable 
network is totally minimal and metrizable.
\qed
\end{corollary}

A topological group is {\em $h$-complete} if all its continuous
homomorphic images are complete (see \cite{DikTon}). 

Using Theorem~\ref{thm:open}, we obtain the following strengthening of the
Guran's Embedding Theorem for sequentially $h$-complete groups:

\begin{theorem} \label{thm:struct}
Every $\omega$-precompact sequentially $h$-complete group $G$ densely 
embeds into the (projective) limit of its metrizable quotients. In 
particular, if $G$ is $h$-complete, then it is equal to the limit of its 
metrizable quotients.
\end{theorem}

\begin{proof}
By Fact~\ref{fact:Guran}, since $G$ is $\omega$-precompact, it embeds into 
$\Sigma = \prod\limits_{\alpha \in I} \Sigma_\alpha$, a product of second 
countable groups. Denote by  $\pi_\alpha: G \rightarrow 
\Sigma_\alpha$ the restriction of the canonical projections to $G$; 
without loss of generality, we may assume that the $\pi_\alpha$ are onto. 
Since $G$ is also sequentially $h$-complete,  the $\pi_\alpha$ are open
(by Theorem~\ref{thm:open}). Thus, one has $G/N_\alpha \cong \Sigma_\alpha$, 
where $N_\alpha = \ker \pi_\alpha$. We may also assume that {\em all}
metrizable quotients of $G$ appear in the product constituting $\Sigma$, 
because by adding factors we do not ruin the embedding.

Let $\iota$ be the embedding of $G$ into the product of its metrizable 
quotients, and put $L = \lim\limits_{\genfrac{}{}{0pt}1\longleftarrow 
{\alpha \in I}}G/N_\alpha$. The image of $\iota$ is obviously contained in 
$L$. (We note that if $G/N_\alpha$ and $G/N_\beta$ are metrizable 
quotients, then $G/N_\alpha N_\beta$ and $G/N_\alpha \cap N_\beta$ are 
also metrizable quotients; the latter one is metrizable, because the 
continuous homomorphism 
$G \rightarrow G/N_\alpha \times G/N_\beta$ is open onto its image, 
as the codomain is metrizable.)  In order to show density, let 
$x= (x_\alpha N_\alpha)_{\alpha \in I} \in L$ and let 
\[
U = U_{\alpha_1} \times \cdots \times U_{\alpha_k} \times
\prod\limits_{\alpha \neq \alpha_i} G/N_\alpha
\]
be a neighborhood of $x$. By the consideration above, the quotient 
$G/ \bigcap\limits_{i=1}^k N_{\alpha_i}$ is metrizable, 
so $\bigcap\limits_{i=1}^k N_{\alpha_i}=N_\gamma$ for some $\gamma \in I$. 
Thus, $\pi_{\alpha_i}(x_\gamma) = x_{\alpha_i}N_{\alpha_i}$, therefore
$\iota(G)$ intersects $U$, and hence $\iota(G)$ is dense in $L$.
\end{proof}

A topological group $G$ is {\em maximally almost-periodic} (or briefly, 
{\em MAP})  if it admits a continuous monomorphism $m: G \rightarrow K$  
into a compact  group $K$, or equivalently, if the  finite-dimensional 
unitary representations of  $G$ separate points.

The Theorem below is a far-reaching generalization of the main result 
of \cite{cikk4}:

\begin{theorem} \label{thm:main}
Let $G$ be a complete topological group. Then the following assertions are 
equivalent:

\begin{list}{{\rm (\roman{enumi})}}
{\usecounter{enumi}\setlength{\labelwidth}{25pt}}

\item
the closed normal subgroups of closed separable subgroups of $G$ are 
$h$-complete and MAP;

\item 
every closed separable subgroup $H$ of $G$ is sequentially $h$-complete, 
and its metrizable quotients $H/N$ are MAP;

\item
$G$ is compact.
\end{list}
\end{theorem}

The following easy consequence of a result by Dikranjan and Tka{\v{c}}enko 
plays a very important role in proving  Theorem~\ref{thm:main}:

\begin{fact} \label{fact:DikTka}
$G$ is precompact if and only if every closed separable subgroup of 
$G$ is precompact. {\rm (Theorem 3.5 in \cite{DikTka}.)}
\end{fact}

\begin{proof}
(i) $\Rightarrow$ (ii): If $H/N$ is a quotient described in (ii), 
then it  is clearly $h$-complete, as the homomorphic image of $H$, which 
is assumed  to be $h$-complete in (i). Let $m: H \rightarrow K$ 
be a continuous  injective homomorphism into a compact group $K$. Since 
$H$ is $h$-complete, $m(H)$ is closed in $K$, so we may assume that $m$ 
is onto. The subgroup $m(N)$ is normal in $K$, because $m$ is bijective, 
and it is closed  because $N$ is $h$-complete. Therefore,  
$\bar m: H/N \rightarrow K/m(N)$ is a continuous injective homomorphism, 
showing that $H/N$ is MAP.

(ii) $\Rightarrow$ (iii): 
Since $G$ is complete, in order to show that $G$ is compact, we show 
that it is also precompact. By Fact~\ref{fact:DikTka}, it suffices to 
show that every closed separable subgroup of $G$ is precompact.

Let $H$ be a closed separable subgroup of $G$. The group $H$ is
$\omega$-precompact (because it is separable), and by (ii) $H$ is
sequentially $h$-complete. Applying Theorem~\ref{thm:struct}, $H$ densely
embeds into the the product $P$ of its metrizable quotient. In order to
show that $H$ is precompact, we show that each factor of the product $P$
is precompact.

Let $Q=H/N$ be a metrizable quotient of $H$; since $H$ is sequentially 
$h$-complete, so is $Q$. The group $Q$ is second countable, because (being 
the continuous image of $H$) it is separable.  Thus, by  
Corollary~\ref{cor:totmin}, $Q$ is totally minimal, because it is 
sequentially  $h$-complete. According to (ii), $Q$ is MAP, and together 
with minimality this implies that $Q$ is precompact. 
\end{proof}

\begin{remark*}
The implication (i) $\Rightarrow$ (iii) can also be proved directly, 
without applying Theorem~\ref{thm:struct}, by using
Fact~\ref{fact:DikTka} and \cite[3.4]{DikUsp}.
\end{remark*}

Since subgroups and  continuous homomorphic images of $c$-compact groups
are $c$-compact again, and are in particular $h$-complete, we obtain:

\begin{corollary} \label{cor:compact}
If $G$ is $c$-compact and MAP, then $G$ is compact.
\end{corollary}

\begin{proof}
If $G$ is MAP, then in particular all its subgroups are so, and all 
its closed subgroups are $c$-compact, thus $h$-complete.
Hence Theorem~\ref{thm:main} applies.
\end{proof}

A group is {\em minimally almost periodic} (or briefly, {\em m.a.p.})
\label{page:map} if it has no non-trivial finite-dimensional 
unitary representations.

\begin{corollary} 
Every $c$-compact group $G$ has a maximal compact quotient $G/M$, where 
$M$ is a closed characteristic m.a.p. subgroup of $G$.
\end{corollary}

\begin{proof}
Let $G$ be a $c$-compact group, and let $M=n(G)$, the von Neumann radical 
of $G$ (the intersection of the kernels of all the finite-dimensional 
unitary representations of $G$). By its definition, $G/M$ is the maximal 
MAP quotient of $G$, and according to Corollary~\ref{cor:compact} this is 
the same as the maximal compact quotient of $G$, because each quotient of 
$G$ is $c$-compact.

The subgroup $M$ is also $c$-compact, so by what we have proved so far, it 
has a  maximal compact quotient $M/L$, where $L$ is characteristic in $M$, 
thus $L$ is normal in $G$. By the Third Isomorphism Theorem, 
$G/M \cong (G/L)/(M/L)$, and since both $M/L$ and $G/M$ are compact, by 
the three space property of compactness in topological groups, the 
quotient $G/L$ is compact too. Thus, $M \subseteq L$, and therefore $M=L$. 
Hence, $M$ is minimally almost periodic.
\end{proof}

We conclude with a reduction theorem. If $G$ is $c$-compact, then so is
every closed subgroup of $G$. Thus, each closed subgroup $H$ and separable
metrizable quotient $H/N$ mentioned in (ii) of Theorem~\ref{thm:main} is
also $c$-compact. Therefore, the equivalence of (ii) and (iii) in
Theorem~\ref{thm:main} yields:

\begin{theorem} \label{thm:reduction}
The following statements are equivalent:

\begin{list}{{\rm (\roman{enumi})}}
{\usecounter{enumi}\setlength{\labelwidth}{25pt}\setlength{\topsep}{2pt}
\setlength{\itemsep}{0pt}}

\item
every $c$-compact group is compact;

\item
every second countable $c$-compact group is MAP (and thus compact);

\item
every second countable $c$-compact m.a.p. group is trivial.
\qed

\end{list}
\end{theorem}

It is not known whether (ii) is true, but by Corollary~\ref{cor:totmin} 
every second-countable $c$-compact group is totally minimal, a fact that
may help in proving or giving a counterexample to (ii).

\section*{Acknowledgements}

I am deeply indebted to Prof. Walter Tholen, my PhD thesis supervisor,
for giving me incredible academic and moral help. 

I am thankful to Prof. Mikhail Tka{\v{c}}enko for sending me his papers
\cite{TkaTop1} and \cite{TkaTop2} by e-mail.

I am grateful for the constructive comments of the anonymous referee that
led to an improved presentation of the results in this paper.

Last but not least, special thanks to my students, without whom it would
not be worth it.

\bibliography{hseq}

\end{document}